\documentclass{amsart}

\usepackage{amsmath}
\usepackage{amssymb}
\usepackage{epsfig}
\usepackage{dsfont}
\usepackage[english]{babel}
\usepackage{amsfonts}

\newcommand{\utilde}[1]{\underset{\!\widetilde{}}{#1}}

\title{Classical Mechanics of Minimal Tori in $\mathbb S^3$}

\author{Joakim Arnlind}
\address[Joakim Arnlind]{Dept. of Math.\\
Link\"oping University
581 83 Link\"oping\\
Sweden}
\email{joakim.arnlind@liu.se}

\author{Jaigyoung Choe}
\address[Jaigyoung Choe]{Korea Institute for Advanced Study\\
207-43 Cheongryangri 2-dong\\
Dongdaemun-gu, Seoul 130-722, Korea}
\email{choe@kias.re.kr}

\author{Jens Hoppe}
\address[Jens Hoppe]{Korean Institute for Advanced Study, 
Royal Institute of Technology,
Sogang University
}

\begin{document}

\begin{abstract}
  We formulate a class of minimal tori in $\mathbb S^3$ in terms of
  classical mechanics, reveal a curious property of the Clifford
  torus, and note that the question of periodicity can be made more
  explicit in a simple way.
\end{abstract}

\maketitle

\section{Introduction}
The Clifford torus $\mathbb S^1(1/\sqrt{2})\times\mathbb
S^1(1/\sqrt{2})$ is geometrically the simplest torus. It is a flat
square torus. When embedded in the three-dimensional unit sphere
$\mathbb S^3$, it becomes a minimal surface. It divides $\mathbb S^3$
into two congruent solid tori. Moreover, it is doubly ruled, i.e. it
can be foliated by two orthogonal families of great circles. The
Clifford torus is the unique algebraic minimal surface of degree 2 and
is characterized even locally as the only (non-totally-geodesic)
minimal surface of contant curvature in $\mathbb S^3$. A long-standing
conjecture by Lawson \cite{L2} that the Clifford torus is the only
embedded minimal torus in $\mathbb S^3$ was recently solved
affirmatively by S. Brendle \cite{B1}. Also, more recently the
Clifford torus has been used as a building block for constructing
infinitely many compact embedded minimal surfaces in $\mathbb S^3$
\cite{KY}, \cite{CS}. For a good recent overview concerning minimal
surfaces in $\mathbb S^3$ see \cite{B2} (and for constant mean
curvature tori see \cite{AL} and references therein).

In this note, we will report on two previously unnoticed aspects of
minimal tori in $\mathbb S^3$: as deformations of the Clifford torus,
and a relation to classical mechanics.

\section{Hamiltonian description of surfaces}
Originally starting out to find minimal surfaces of higher genus via stationary points of the functional
\begin{eqnarray}
S[w]:=\int\sqrt{\eta^{KL}\partial_{K}w\partial_{L}w}\sqrt{\eta}d^{D}u
\end{eqnarray}
i.e. solutions of (c.p. \cite{BH}, e.g. )
\begin{eqnarray}
(\eta^{IJ}\eta^{KL}-\eta^{IK}\eta^{JL}) \partial_{I}w\partial_{J}w (\partial^{2}_{KL}-\Gamma^{M}_{KL}\partial_{M}w)
=0,
\end{eqnarray}
we found a class of minimal tori of the form
\begin{eqnarray}
\vec{x}(\varphi^{1},\varphi^2)=\left(\begin{array}{cc}
                                            \cos\theta(\varphi^{1},\varphi^2)\cos(\varphi^{1}) \\
                                            \cos\theta(\varphi^{1},\varphi^2)\sin(\varphi^{1}) \\
                                            \sin\theta(\varphi^{1},\varphi^2)\cos(\varphi^{2}) \\
                                            \sin\theta(\varphi^{1},\varphi^2)\sin(\varphi^{2})
                                          \end{array}\right) \in \mathbb S^3 \subset \mathbb{R}^{4},
\end{eqnarray}
where $\theta$ (when depending on $\varphi^{1}$ and $\varphi^{2}$ only via the combination $k\varphi^{1}+l\varphi^{2}=:t$) can be determined as the solution of a classical mechanics problem,
\begin{eqnarray}
\dot{\theta}^{2}+\frac{c^{2}s^{2}}{k^{2}s^{2}+l^{2}c^{2}}\left(1-\frac{c^{2}s^{2}}{E^{2}}\right)=0,
\end{eqnarray}
i.e. a zero energy solution of a point-particle (having "position" $\theta$ at ``time'' $t$) moving in a potential $V_E(\theta)$. ($c:=\cos\theta$, $s:=\sin\theta$)

In order to see that solutions of (4) give minimal surfaces of the form (3), let us calculate the first and second fundamental forms corresponding to (3):
\begin{eqnarray}
(g_{ab})&=& \left(\begin{array}{cc}
              c^{2}+(\partial_{1}{\theta})^{2} & \partial_{1}{\theta}\partial_{2}{\theta} \\
              \partial_{1}{\theta}\partial_{2}{\theta} & s^{2}+(\partial_{2}{\theta})^{2} \\
            \end{array}\right)\\
(h_{ab})&=& \frac{1}{|\overrightarrow{m}|}\left(\begin{array}{cc}
              s^{2}c+sc\theta_{11}+2S^{2}\theta_{1}^{2} & sc\theta_{12}+(s^{2}-c^{2})\theta_{1}\theta_{2} \\
              sc\theta_{12}+(s^{2}-c^{2})\theta_{1}\theta_{2} & -s^{2}c^{2}+sc\theta_{22}-2c^{2}\theta^{2}_{2} \\
            \end{array}\right)
\end{eqnarray}
with
\begin{eqnarray}
\overrightarrow{m}:=sc\left(
                        \begin{array}{c}
                          \begin{array}{c}
                              -sc_{1} \\
                              -ss_{1}
                            \end{array}
                           \\
                          \begin{array}{c}
                             cc_{2} \\
                             cs_{2}
                           \end{array}
                           \\
                        \end{array}
                      \right)
                    -\left(\begin{array}{c}

                                                              -ss_{1}\theta_{1} \\
                                                              sc_{1}\theta_{1} \\
                                                              -cs_{2}\theta_{2} \\
                                                              cc_{2}\theta_{2} \\

                                             \end{array}\right),\,\,\,c_i:=\cos\varphi^i,\,s_i:=\sin\varphi^i
\end{eqnarray}
being orthogonal to $\partial_{1}\vec{x}$, $\partial_{2}\vec{x}$, and $\vec{x}$;
hence one finds that (3) has zero mean curvature in $\mathbb S^3$ if and only if
\begin{eqnarray}
\begin{array}{c}
  \left(s^{2}+\theta_{2}^{2} \right)\left(s^{2}c^{2}+sc\theta_{11}+2s^{2}\theta_{1}^{2} \right) + \left(c^{2}+\theta_{1}^{2} \right)\left(-s^{2}c^{2}+sc\theta_{22}-2c^{2}\theta_{2}^{2} \right) \\
 \quad- 2\theta_{1}\theta_{2}\left( sc\theta_{12}+(s^{2}-c^{2})\theta_{1}\theta_{2}\right)=0,
\end{array}
\end{eqnarray}
where $\theta_{a}:=\partial_{a}\theta$.

Note that, equivalently, one could have obtained (8) by varying
\begin{eqnarray}
S[\theta]:= \int\sqrt{g}d\varphi^{1}d\varphi^{2}=\int\sqrt{c^{2}s^{2}+s^{2}\theta_{1}^{2}+c^{2}\theta^{2}_{2}}d\varphi^{1}d\varphi^{2} .
\end{eqnarray}
For $\theta =\theta(k\varphi^{1}+l\varphi^{2})=\theta(t)$ one gets, from the Lagrangian (the `$-$' sign put in for later convenience)
\begin{eqnarray}
L= -\sqrt{c^{2}s^{2}+(k^{2}s^{2}+l^{2}c^{2})\dot{\theta}^{2}}
\end{eqnarray}
as well as directly from (8), the second order equation
\begin{eqnarray}
\begin{array}{c}
  sc\ddot{\theta}(k^{2}s^{2}+l^{2}s^{2})+ \dot{\theta}^{2}\left[(l^{2}-k^{2})s^{2}c^{2}+2s^{4}k^{2}-2c^{4}l^{2}\right] \\
 \quad+s^{2}c^{2}(s^{2}-c^{2})=0 ,
\end{array}
\end{eqnarray}
which can also be shown to follow from (4) by simply differentiating. A possible way to systematically derive (4) from (11) is provided by a standard Legendre-transformation, obtaining from (10) a Hamiltonian
\begin{equation}
H:=\frac{\partial L}{\partial \dot{\theta}}\dot{\theta}-L=\cdots=\frac{c^{2}s^{2}}{\sqrt{c^{2}s^{2}+(k^{2}s^{2}+l^{2}c^{2})\dot{\theta}^{2}}}=E.
\end{equation}
Expressing H in terms of $\theta$ and the canonical momentum
\begin{eqnarray}
\pi:=\frac{\partial L}{\partial \dot{\theta}}=-\frac{(k^{2}s^{2}+l^{2}s^{2})\dot{\theta}}{\sqrt{c^{2}s^{2}+(k^{2}s^{2}+l^{2}c^{2})\dot{\theta}^{2}}}
\end{eqnarray}
one obtains (replacing $|cs|$ by $cs$, justifyable via $H = const$)
\begin{eqnarray}
H= \frac{1}{2}\sin(2\theta)\sqrt{1-\frac{\pi^{2}}{ k^{2}s^2+l^2c^2}}=H[\theta,\pi] ;
\end{eqnarray}
one can check that the first-order equations
\begin{equation}
  \dot{\theta}=\frac{\partial H}{\partial\pi}= -\frac{1}{2}\sin 2\theta \frac{\frac{\pi}{ k^{2}s^2+l^2c^2}}{\sqrt{1-\frac{\pi^{2}}{ k^{2}s^2+l^2c^2}}}\,;\quad
\dot{\pi} = -\frac{\partial H}{\partial \theta}
\end{equation}
reproduce (11).

\section{$k=l$: Clifford torus in disguise}
Whereas the axially symmetric $k=0$ (or $l=0$) case, except for our mechanical interpretation, is well known (S.Brendle \cite{B2} quotes R. Kusner, when discussing an elliptic integral solution of the form (3) for $l=1, k=0$), the above mentioned $k=l$ case can actually be integrated in terms of elementary functions, leading to the curious fact that the Clifford-Torus can be viewed as a "non-trivial" graph over itself (in infinitely many different ways).
Separating variables, and letting
$$\alpha\left(\varphi:= \varphi^{1}+\varphi^{2}\right)=2\theta\left(k\left(\varphi^{1}+\varphi^{2}\right)\right),$$
$a:= \frac{1}{2E}=\cosh\gamma >1$, one deduces from (4) that
%\begin{eqnarray} \nonumber \\
%\begin{array}{c}
\begin{equation}  \frac{t-t_{0}}{k}=\varphi-\varphi_{0}=\pm \int\frac{d\alpha}{\sin\alpha\sqrt{a^{2}\sin^{2}\alpha-1}}
\end{equation}
\begin{equation*}
=\mp\arctan\frac{\cos\alpha}{\sqrt{a^{2}\sin^{2}\alpha-1}}\,\,\,\,\,\,\,\,\,\,
%\end{array}
%\end{eqnarray}
\end{equation*}
i.e. that $\alpha(\varphi)$ is given via
\begin{eqnarray}
\frac{\mp\cos\alpha}{\sqrt{a^{2}\sin^{2}\alpha-1}}=\tan(\varphi-\varphi_{0})=:T
\end{eqnarray}
as a solution of
\begin{eqnarray}
(\alpha')^{2}=\sin^2\alpha \left(a^{2}\sin^{2}\alpha-1\right).
\end{eqnarray}
Using this, respectively
\begin{equation}
  \sin\alpha=\sqrt{\frac{1+T^{2}}{1+a^{2}T^{2}}}=\frac{1}{\sqrt{\cos^{2}(\varphi-\varphi_{0})+\cosh^{2}\gamma\sin^{2}(\varphi-\varphi_{0})}}
\end{equation}
\begin{equation*}
 \cos\alpha=\frac{-\sin(\varphi-\varphi_{0})\sinh\gamma}{\sqrt{\cos^{2}(\varphi-\varphi_{0})+\cosh^{2}\gamma\sin^{2}(\varphi-\varphi_{0})}} ,\,\,\,\,\,\,\,\,\,\,\,\,\,\,\,\,\,\,\,\,\,\,\,\,\,\,\,
\end{equation*}
one can check that $\alpha(\varphi)$, given as in (19), does satisfy the second order equation (cp. (11)) that is equivalent to the vanishing of the mean curvature,
\begin{eqnarray}
(\sin \alpha) \alpha''-2\alpha'^{2}\cos\alpha=(\cos\alpha)(\sin\alpha)^{2},
\end{eqnarray}
whereby it is useful to note that
$$
\alpha'\cot\alpha= -\frac{T(a^{2}-1)}{1+a^{2}T^{2}}
$$
\begin{eqnarray}
\alpha''=\frac{-\sin(\varphi-\varphi^{0})\sinh\gamma(\cosh^{2}\gamma+\sinh^{2}
\gamma\cos^{2}(\varphi-\varphi^{0}))}{(\cos^{2}(\varphi-\varphi^{0})+\cosh^{2}
\gamma\sin^{2}(\varphi-\varphi^{0}))^{2}}.
\end{eqnarray}
The derived solution(s) read(s)
\begin{eqnarray}
\vec{x}\left(\varphi^{1}\varphi^{2}\right)=\frac{1}{2}\left(
  \begin{array}{c}
    \sqrt{1-\frac{e \sin(\varphi-\varphi^{0})}{\sqrt{1+e^{2}\sin^{2}(\varphi-\varphi^{0})}}}\begin{array}{c}
    \cos\varphi^{1} \\
    \sin\varphi^{1}
    \end{array}\\
    \sqrt{1+\frac{e \sin(\varphi-\varphi^{0})}{\sqrt{1+e^{2}\sin^{2}(\varphi-\varphi^{0})}}}\begin{array}{c}
    \cos\varphi^{2} \\
    \sin\varphi^{2}
    \end{array} \\
  \end{array}
\right)
\end{eqnarray}
with $e:= \sinh\gamma \in \mathbb{R}$ an arbitrary constant.

One easily sees that, for each value of $e \in \mathbb{R}$, the solution (22) defines an embedded minimal torus in $S^{3}$, which at first sight is rather puzzling, as every embedded minimal torus in $\mathbb S^3$ must be congruent to the Clifford torus \cite{B1}. Let us give three direct proofs for the above concrete case, (22).

Firstly, for constant $\varphi$ both $x_{3}$ and $x_{4}$ can be expressed as a linear combination of $x_{1}$ and $x_{2}$, defining 2 hyper-planes in $\mathbb{R}^{4}$, respectively. Since the intersection of $\mathbb S^3$ with two hyper-planes containing the origin of $\mathbb R^4$  gives a great circle, the minimal torus is ruled. It is known \cite{L1} that the Clifford torus is the only embedded surface among all infinitely many ruled minimal surfaces in $\mathbb S^3$. Hence (22) describes a  Clifford torus.

Secondly, calculating the determinants of (6) and (5) for $\theta = \theta\left(k\varphi^{1}+l\varphi^{2}\right)$, one finds
\begin{equation}
  h=\frac{-\dot{\theta}^{4}-2s^{2}c^{2}\dot{\theta}^{2}-s^{4}c^{4}}{\dot{\theta}^{2}+s^{2}c^{2}}=-\left(\dot{\theta}^{2}+s^{2}c^{2}\right)
\end{equation}
\begin{equation*}
  g=\dot{\theta}^{2}+s^{2}c^{2},\qquad \frac{h}{g}=-1 ,\qquad\qquad\qquad\,\,\,
\end{equation*}
hence for the intrinsic Gaussian curvature
\begin{eqnarray}
\begin{array}{c}
  R= (Tr W)^{2}-\left(Tr W^{2}\right)+(Tr \tilde{W})^{2}-\left(Tr \tilde{W}^{2}\right) \\
  =0 - 2 + 4 -2=0\qquad\qquad\qquad\qquad\,\,\,\,\,\,\,\,\,\,\,
\end{array}
\end{eqnarray}
$$\left(\tilde{h}_{ab}:= \vec{x}\cdot\partial _{ab}^{2}\vec{x}=-\partial _{a}\vec{x}\cdot\partial _{b}\vec{x}=-g_{ab}\Rightarrow \tilde{W}^{a}_{b}:=g^{ac}\tilde{h}_{cb}=-\delta^{a}_{b}\right),$$ proving that the minimal tori (22) are flat (hence the Clifford torus).

Thirdly, one can construct an explicit isometry to
\begin{eqnarray}
\vec{\utilde{x}}(\tilde{\varphi}^{1},\tilde{\varphi}^2)=\frac{1}{\sqrt{2}}\left(\begin{array}{c}
                                                                           \cos(\tilde{\varphi}^{1}) \\
                                                                           \sin(\tilde{\varphi}^{1}) \\
                                                                          \cos(\tilde{\varphi}^{2}) \\
                                                                          \sin(\tilde{\varphi}^{2})
                                                                         \end{array}\right)
\end{eqnarray}
by finding a re-parametrization  $\varphi^{1}\varphi^{2}\rightarrow\tilde{\varphi}^{1} \tilde{\varphi}^2$ such that $(k=l=1, \varphi^{0}=0)$
\begin{equation}
  2\left(g_{ab} \right)=  \mathds{1} + \left(
                              \begin{array}{cc}
                                \cos\alpha & 0 \\
                                0 & -\cos\alpha \\
                              \end{array}
                            \right)
  +\frac{1}{2}\sin^{2}\alpha\left(a^{2}\sin^{2}-1\right)\left(
                                                          \begin{array}{cc}
                                                            1 & 1 \\
                                                            1 & 1 \\
                                                          \end{array}
                                                        \right)
 \end{equation}
\begin{equation*}
  \qquad\qquad\quad\,\,\,\, = \mathds{1} -\frac{e\sin\varphi}{\sqrt{1+e^{2}\sin^{2}\varphi}}\left(
                                                       \begin{array}{cc}
                                                         1 & 0 \\
                                                         0 & -1 \\
                                                       \end{array}
                                                     \right) + \frac{1}{2}\frac{e^{2}\sin^{2}\varphi}{(1+e^{2}\sin^{2}\varphi)^{2}}
                                                     \left(
                                                       \begin{array}{cc}
                                                         1 & 1 \\
                                                         1 & 1 \\
                                                       \end{array}
                                                     \right)
\end{equation*}
\begin{equation*}
  \quad= J^{T}(2\tilde{g}_{..})J=J^{T}J;\qquad J^{\tilde{a}}_{a}=\frac{\partial \tilde{\varphi}^{\tilde{a}}}{\partial\varphi^{a}}.\,\,\,\,\,\,\,\,\,\,\,\,\,\,\,
  \,\,\,\,\,\,\,\,\,\,\,\,\,\,\,\,\,\,\,\,\,\,\,\,\,\,\,\,\,\,\,\,
\end{equation*}
The Ansatz
\begin{eqnarray}
\tilde{\varphi}^{1}=\varphi^{1}+\int^{\varphi}u,\qquad \tilde{\varphi}^{2}=\varphi^{2}+\int^{\varphi}v
\end{eqnarray}
gives
\begin{eqnarray}
\begin{array}{c}
  J=\left(
    \begin{array}{cc}
      1+u & u \\
      v & 1+v \\
    \end{array}
  \right), \qquad \qquad\qquad\qquad\qquad\,\,\\
  J^{T}J= \mathds{1} + \left(
                              \begin{array}{cc}
                                2u & u+v \\
                                u+v & 2v \\
                              \end{array}
                            \right)+ (u^{2}+v^{2}) \left(
                                                     \begin{array}{cc}
                                                       1 & 1 \\
                                                       1 & 1 \\
                                                     \end{array}
                                                   \right) ,
\end{array}
\end{eqnarray}
so that (26) will be satisfied when choosing
\begin{equation}
  u=\frac{1}{2}\cos\alpha(\varphi)+ w(\varphi) \,\,\,\,\,
\end{equation}
\begin{equation*}
  v=-\frac{1}{2}\cos\alpha(\varphi)+ w(\varphi)
\end{equation*}

\begin{equation*}
  u^{2}+v^{2}+2w = w^{2} + 2w + \frac{1}{4}(\cos\alpha)^{2}=\frac{1}{2}\frac{c^{2}\sin^{2}\varphi}{(1+c^{2}\sin^{2}\varphi)^{2}}\qquad \end{equation*}
  \begin{equation*}
  w(\varphi)=-1+\frac{\sqrt{1+\frac{3}{4}e^{4}s^{4}+\frac{9}{4}e^{2}s^{2}}}{\left(1+e^{2}s^{2}\right)}        (s=\sin\varphi).
\end{equation*}

\section{$k\neq l$}
Analogously one could explicitly construct isothermal coordinates also for arbitrary $k$ and $l$ via $(s=\sin\theta,c=\cos\theta)$
\begin{eqnarray}
\begin{array}{c}
  \left(g_{ab}\right)=\left(
                      \begin{array}{cc}
                        \cos2\theta & 0 \\
                        0 & \sin2\theta \\
                      \end{array}
                    \right)+\left(
                              \begin{array}{cc}
                                k^{2} & kl \\
                                kl & l^{2} \\
                              \end{array}
                            \right)\left(\frac{c^{2}s^{2}}{E^{2}}-1\right)\left(\frac{c^{2}s^{2}}{k^{2}s^{2}+l^{2}c^{2}}\right) \\
  = \rho^{2}J^{T}J, \qquad J=\left(
                                     \begin{array}{cc}
                                       1+ku & lu \\
                                       k v & 1+l v  \\
                                     \end{array}
                                   \right) .\qquad\qquad\qquad
\end{array}
\end{eqnarray}
With
$$Y(\theta):=\left(\frac{c^{2}s^{2}}{E^{2}}-1\right)\left(\frac{c^{2}s^{2}}{k^{2}s^{2}+l^{2}c^{2}}\right) = -V_E(\theta)$$
one gets
\begin{eqnarray}
\begin{array}{c}
  c^{2}+k^{2}Y=\left(1 + k^{2}(u^{2}+v^{2})+2ku\right)\rho^{2} \\
  s^{2}+l^{2}Y=\left(1 + l^{2}(u^{2}+v^{2})+2lv\right)\rho^{2}\,\, \\
  \qquad\, klY=\left(lu + k v + kl(u^{2}+v^{2})\right)\rho^{2} .
\end{array}
\end{eqnarray}
Substracting $2kl(klY)$ from $l^{2}(c^{2}+k^{2}Y) + k^{2}(s^{2}+l^{2}Y) $ and dividing by $k^{2} + l^{2}$ one obtains
\begin{eqnarray}
\rho^{2}=\frac{k^{2}s^{2}+l^{2}c^{2}}{k^{2} + l^{2}}
\end{eqnarray}
while $l^{2}(c^{2}+k^{2}Y) - k^{2}(s^{2}+l^{2}Y) $ yields an (inhomogeous) linear! relation between $u$ and $v$, which (when substituted into any of the 3 equations in (31)) gives a simple quadratic equation for $u$ (or $v$).\\

In order to study the periodicity properties (for general $k$ and $l$), it is convenient to still write the solution into the form (22), but with $\varphi$ replaced by $\tilde{\varphi}$, which is defined by the absorption of the factor $k^{2}s^{2}+l^{2}c^{2}$, i.e.
\begin{eqnarray}
\frac{\partial\varphi}{\partial\tilde{\varphi}}=\sqrt{k^{2}\sin^{2}\tilde{\theta} + l^{2}\cos^{2}\tilde{\theta}}\qquad\qquad\qquad\qquad\,
\end{eqnarray}
$$=\sqrt{\frac{k^{2} + l^{2}}{2}}\sqrt{1+\frac{(k^{2} - l^{2})}{k^{2} + l^{2}}\frac{e\sin\tilde{\varphi}}{\sqrt{1+e^{2}\sin^{2}\tilde{\varphi}}}}.$$
In the case discussed by Brendle ($k=0, l=1$) one would get
\begin{eqnarray}
\varphi^{2}=\frac{1}{\sqrt{2}}\int^{\tilde{\varphi}}\sqrt{1 - \frac{e\sin u}{\sqrt{1 + e^{2}\sin^{2}u}}} du,
\end{eqnarray}
which is a somewhat simpler elliptic integral than the formula for the period in \cite{B2}.
As for $e= 0$,
\begin{eqnarray}
\delta \varphi^{2}= \frac{1}{\sqrt{2}} \int_{\tilde{\varphi}}^{\tilde{\varphi}+2\pi} = \sqrt{2}\pi = \delta_0\varphi^{2},
\end{eqnarray}
while for $e\rightarrow +\infty$,
\begin{eqnarray*}
\delta \varphi^{2}= \frac{1}{\sqrt{2}} \int_{\pi}^{2\pi}\sqrt{2}= \pi
\end{eqnarray*}
it is, because of the continuity in $e$, clear that there will be infinitely many values of $e$ for which $\delta \varphi^{2}$ will be a rational multiple of $\pi$, $(p/q)\,\pi$. Alternatively, with $\sin^2\theta(\varphi^{2}) = r^{2}(\varphi^{2})= v(\varphi^{2})$ one has
\begin{equation}
  \varphi^{2} - \varphi_{0}^{2}= \frac{1}{2}\frac{1}{\sqrt{1- v_{-}}}\pi\left(\arcsin\sqrt{\frac{1 - v}{1 - v_{+}}},1-v_{+},\sqrt{\frac{1- v_{+}}{1- v_{-}}}\right)
\end{equation}
\begin{equation*}
  1 > v \geq v_{+} > v_{-} >0, v_{\pm} = \frac{1}{2}\pm\sqrt{\frac{1}{4} - E^{2}},
\end{equation*}
an elliptic integral of the third kind.

\section*{Acknowledgment}

This work was supported by the Swedish Research Council,
Postech SRC-GAIA, and a Sogang University
Research Grant (2012).


\begin{thebibliography}{11111}
\setlength{\baselineskip}{0.8\baselineskip}
\bibitem[AL]{AL} B. Andrews and H. Li, {\it Embedded constant mean curvature tori in the three-sphere}, arXiv: math.DG/1204.5007v3
\bibitem[BH]{BH} M. Bordemann and J. Hoppe, {\it The dynamics of relativistic membranes. II. Nonlinear waves and covariantly reduced membrane equations}, Phys. Lett. B {\bf 325} (1994), no. 3-4, 359--365.
\bibitem[B1]{B1} S. Brendle, {\it Embedded minimal tori in $\mathbb S^3$ and the Lawson conjecture}, to appear in Acta Math.
\bibitem[B2]{B2} S. Brendle, {\it Minimal surfaces in $\mathbb S^3$ $-$ A survey of recent results},
Bull. Amer. Math. Soc. (2) {\bf } (2013), 133--171.
\bibitem[CS]{CS} J. Choe and M. Soret, {\it New minimal surfaces in $\mathbb S^3$ desingularizing the Clifford tori}, preprint.
\bibitem[KY]{KY} N. Kapouleas and S.-D. Yang, {\it Minimal surfaces in the three-sphere by doubling
the Clifford torus}, Amer. J. Math. (2) {\bf 132} (2010), 257--295.
\bibitem[L1]{L1} H.B. Lawson, Jr., {\it Complete minimal surfaces in $\mathbb{S}^3$}, Ann. Math.
     {\bf 92} (1970), 335--374.
\bibitem[L2]{L2} H.B. Lawson, Jr., {\it The unknottedness of minimal embeddings}, Invent. Math. {\bf 11}  (1970), 183--187.
\end{thebibliography}
\end{document}